\documentclass[aop,citesort,MSNbibl,dvips]{arximspdf}

\doi{10.1214/10-AOP529}
\volume{38}
\issue{5}
\pubyear{2010}
\firstpage{1986}
\lastpage{2008}

\makeatletter

\newtheorem{theorem}{Theorem}
\newtheorem{lemma}[theorem]{Lemma}
\newtheorem{proposition}[theorem]{Proposition}
\newproclaim{remark}[theorem]{Remark}
\newtheorem{corollary}[theorem]{Corollary}
\newproclaim{definition}{Definition}

\makeatother

\begin{document}
\begin{frontmatter}

\title{Best constants in Rosenthal-type inequalities and the Kruglov
operator\protect\thanksref{T1}}
\runtitle{Best constants in Rosenthal's type inequalities}

\thankstext{T1}{Supported by the ARC.}

\begin{aug}
\author[A]{\fnms{S. V.} \snm{Astashkin}\ead[label=e1]{astashkn@ssu.samara.ru}} and
\author[B]{\fnms{F. A.} \snm{Sukochev}\corref{}\ead[label=e2]{f.sukochev@unsw.edu.au}}
\runauthor{S. V. Astashkin and F. A. Sukochev}
\affiliation{Samara State University and University of New South Wales}
\address[A]{Department of Mathematics\\
\quad and Mechanics\\
Samara State University\\
443011 Samara, Acad. Pavlov, 1\\
Russian Federation\\
\printead{e1}} 
\address[B]{School of Mathematics and Statistics\\
University of New South Wales\\
Kensington NSW 2052\\
Australia\\
\printead{e2}}
\end{aug}

\received{\smonth{7} \syear{2009}}
\revised{\smonth{12} \syear{2009}}

%
\begin{abstract}
Let $X$ be a symmetric Banach function space on $[0,1]$ with the
Kruglov property, and let $\mathbf{f}=\{f_k\}_{{k=1}}^n$, $n\ge1$ be an
arbitrary sequence of independent random variables in $X$. This paper
presents sharp estimates in the deterministic characterization of the
quantities
\[
\Biggl\| \sum
_{{k=1}}^nf_k \Biggr\|_X, \Biggl\| \Biggl(\sum
_{{k=1}}^n|f_k|^p \Biggr)^{1/p} \Biggr\|_X,\qquad 1\leq p<\infty,
\]
in terms of the sum of disjoint copies of individual terms of
$\mathbf{f}$. Our method is novel and based on the important recent
advances in the study of the Kruglov property through an operator
approach made earlier by the authors. In particular, we discover that
the sharp constants in the characterization above are equivalent to the
norm of the Kruglov operator in $X$.
\end{abstract}

%

\end{frontmatter}

\section{Introduction}\label{sec1}
For an arbitrary sequence $\mathbf{f}:=\{f_k\}_{{k=1}}^n\subset L_1[0,1]$ consider
its disjointification, that is, the function
%
%
\begin{equation}\label{f1}
F(u):=\sum
_{i=1}^n\bar f_i(u)\qquad (u>0),
\end{equation}
where the sequence $\{\bar f_k\}_{k=1}^n$ is a disjointly supported
sequence of equimeasurable copies of the individual elements from
the sequence $\mathbf{f}$ (we always assume here that $[0,1]$ is
equipped with the Lebesgue measure $\lambda$). Denote by $F^*$ the
decreasing rearrangement of $|F|$ (see relevant definitions in
Section \ref{Preliminaries}). Let $X$ be a symmetric Banach function
space on $[0,1]$ for which there exists a universal constant $C_X>0$
such that the inequality
%
%
\begin{equation}\label{Rosenthal}
\Biggl\|\sum_{k=1}^nf_k \Biggr\|_X\leq
C_X\bigl(\bigl\|F^*\chi_{[0,1]}\bigr\|_X+ \bigl\|F^*\chi_{[1,\infty]}\bigr\|_{L_1}\bigr)
\end{equation}
(resp.,
%
%
\begin{equation}\label{meanRosenthal}
\Biggl\|\sum_{k=1}^n f_k\Biggr\|_X\leq C_X \bigl(\bigl\|F^*\chi_{[0,1]}\bigr\|_X+
\bigl\|F^*\chi_{[1,\infty]}\bigr\|_{L_2} \bigr) \Bigg)
\end{equation}
holds for every sequence $\mathbf{f}\subset X$ of independent
random variables (i.r.v.'s) (resp., of mean zero
i.r.v.'s) and for every $n\in\mathbb{N}$. The inequalities above can
be viewed as a fundamental generalization of the famous
Khintchine inequality and, in the case when $X=L_p$
(resp., $X\supset L_p$), $1\le p<\infty$, they may be
found in \cite{Ros} (resp., \cite{JS}). The original
proof of Rosenthal in \cite{Ros}, as well as a subsequent proof
of a more general result by Burkholder in \cite{Bur} yielded
only constants $C_{L_p}$ in (\ref{Rosenthal}) and
(\ref{meanRosenthal}) which grow exponentially in $p$, as $p\to
\infty$. The sharp result that
$C_{L_p}\asymp\frac{p}{\ln(p+1)}$, that is, there are universal
constants $0<\alpha<\beta<\infty$ such that the ratio between
$C_{L_p}$ and $\frac{p}{\ln(p+1)}$ lies in the interval
$[\alpha,\beta]$ for all $p\in[1, \infty)$ was obtained in
\cite{JSZ} (see also subsequent alternative proofs in
\cite{La,Ju}).

The main purpose of this paper is to provide a sharp estimate on
the constant $C_X$ in (\ref{Rosenthal}) in the more general
setting of symmetric spaces $X$. At the same time, we also believe
that our methods shed additional light on the well-studied case
$X=L_p$. Indeed, the methods exploited in \cite{JSZ,La,Ju} have a
distinct $L_p$-flavor and do not appear to extend to other
symmetric function spaces for which (\ref{Rosenthal}) and
(\ref{meanRosenthal}) hold. Our approach here is linked with the
so-called Kruglov property (see the definition in Section
\ref{The Kruglov property}). Consider the special case of
Rosenthal's inequalities (\ref{Rosenthal}) and
(\ref{meanRosenthal}) when i.r.v.'s $f_k$, $k=1,2,\ldots,n$, satisfy
the additional assumption that
%
%
\begin{equation}\label{assumption}
\sum_{k=1}^n\lambda(\{f_k\neq0\} )\leq1,\qquad n\in\mathbb{N}
\end{equation}
[in this case, the right-hand sides of (\ref{Rosenthal}) and
(\ref{meanRosenthal}) become equal]. In this special case, it was
first established by Braverman \cite{Br} that if $X$ is a
symmetric space with the Fatou property (see Section
\ref{spaces} below), then $X$ has the Kruglov property if and only
if (\ref{Rosenthal}) holds. Recently, in \cite{AS1,AS2,AS3} we
have developed a novel approach to the study of spaces with the
Kruglov property that involves defining a positive linear operator
$K\dvtx L_1[0,1]\to L_1[0,1]$ (see details in Section
\ref{The Kruglov operator}) which is bounded in a symmetric function space
$X$ with the Fatou property if and only if $X$ has the Kruglov
property. Furthermore, we have shown in those papers that in this
case the (Kruglov) operator $K$ is bounded in $X$ if and only if
(\ref{Rosenthal}) and (\ref{meanRosenthal}) hold in full
generality. The following key fact is an immediate consequence of
Prohorov's familiar inequality \cite{Pr} (see also the proof of
Theorem 3.5 from \cite{AS1}): if $X$ is a symmetric space, then
for every sequence $\{f_k\}_{k=1}^\infty\subset L_1[0,1]$ of
symmetrically distributed i.r.v.'s satisfying assumption
(\ref{assumption}) we have
\[
\Biggl(\sum_{k=1}^nf_k \Biggr)^*\leq16K(F)^*\qquad (n\in
\mathbb{N}),
\]
where the function $F$ is defined by (\ref{f1})
(the definition of the decreasing rearrangement $f^*$ of a
measurable function $f$ is given in the next section).
This observation has naturally led us to the conjecture that
the best constant $C_X$ in (\ref{Rosenthal}) and
(\ref{meanRosenthal}) should be equivalent to the norm of the
operator $K$ in $X$. We prove this conjecture in Section
\ref{The role} and present computations of the norm $\|K\|$ in
various classes of symmetric spaces in Section \ref{The norm}.
In the case of $L_p$-spaces, $1\le p<\infty$, our results, of
course, yield the same estimates as in \cite{JSZ,La,Ju}. In
the case of symmetric Lorentz and Marcinkiewicz spaces (and
other classes of symmetric spaces in which we are able to
compute the norm of the operator $K$) our results are new and
appear to be unattainable by methods used in \cite{JSZ,La,Ju}.
In the final section of this paper we provide two complements
to Rosenthal's inequality (\ref{Rosenthal}).


\section{Preliminaries}\label{Preliminaries}
\subsection{Symmetric function spaces and interpolation of
operators}\label{spaces}
In this subsection we present some definitions from the theory of
symmetric spaces and interpolation of operators. For more details on
the latter theory we refer to \cite{KPS,LT2}.

We will denote by $S(\Omega,{\mathcal P})$ the linear space of all
measurable finite a.e. functions on a given measure space $(\Omega,
{\mathcal P})$ equipped with the topology of convergence locally in
measure.

Let $I$ denote either $[0,1]$ or $(0,\infty)$ with Lebesgue measure
$\lambda$. If $f\in S(I,\lambda)$ we denote by $f^*$ the decreasing
rearrangement of $f$, that is,
\[
f^*(t)={\inf_{\lambda(A)=t}\sup_{s\in I\setminus A}}|f(s)|.
\]
%

A Banach function space $X$ on $I$ is said to be \textit{symmetric}
if the conditions $f\in X$ and $g^*\le f^*$ imply that $g\in X$ and $\|
g\|_X\le\|f\|_X$.
We will assume always the normalization that $\|\chi_{(0,1)}\|_X=1$, where
$\chi_{A}$ is the characteristic function of the set $A\subset I$.
Let $\varphi_X(t)=\|\chi_{(0,t)}\|_X$ be the \textit{fundamental
function} of $X$.
A symmetric space $X$ is said to have \textit{the Fatou property} if
for every sequence $(f_n)_{n=1}^{\infty}\subset X$ of nonnegative
functions such that $f_n\uparrow f$ a.e. and
$\lim_{n\to\infty}\|f_n\|_X<\infty$ we have $f\in X$ and
$\|f\|_X=\lim_{n\to\infty}\|f_n\|_X$.

Let us recall some classical examples of symmetric spaces on
$[0,1]$.

Let $M(t)$ be an increasing convex function on $[0,\infty)$ such
that $M(0)=0$. By $L_M$ we denote \textit{the Orlicz space} on $[0,1]$
(see, e.g., \cite{KPS,LT2}) endowed with the norm
\[
\|x\|_{L_M}=\inf\biggl\{\lambda>0\dvtx \int_{0}^{1}
M\bigl(|x(t)|/\lambda\bigr)\,dt\leq1\biggr\}.
\]

We suppose that $\psi$ is a positive concave function on $[0,1]$
with $\psi(0+)=0$. \textit{The Lorentz space} $\Lambda(\psi) $ is the
space of all measurable functions $f$ on the interval $[0,1]$ such
that
\[
\Vert f\Vert_{\Lambda(\psi)} =\int_0^1f^*(s)\,d\psi(s)<\infty.
\]
\textit{The Marcinkiewicz space} $M(\psi)$ is the space of all
measurable functions $f$ on the interval $[0,1]$ such that
\[
\Vert f\Vert_{M(\psi)} =\sup_{0<t\le1}\frac{\psi(t)}{t}\int
_0^tf^*(s)\,ds<\infty.
\]
It is easy to check that $\varphi_{\Lambda(\psi)}(t)=
\varphi_{M(\psi)}(t)=\psi(t)$. In this paper, we mainly work with
the case when $\psi(t)=t^{1/p}$, $1\le p<\infty$.

Let $\vec{{X}}=(X_0,X_1)$ be a Banach couple and $X$ be a
Banach space such that $X_0\cap X_1\subseteq X\subseteq X_0+X_1$. We
say that $X$ is \textit{an interpolation space} between $X_0$ and $X_1$
if any bounded linear operator $A\dvtx X_0+X_1\to X_0+X_1$ which maps
$X_i$ boundedly into $X_i$ $(i=0,1)$ also maps $X$ boundedly into
$X$. Then $\|A\|_{X\to X}\le C(\|A\|_{X_0\to X_0},\|A\|_{X_1\to
X_1})$ for some $C\ge1$. If the last inequality holds with $C=1$ we
will refer $X$ to an \textit{1-interpolation space} between $X_0$ and
$X_1$.

In what follows $\operatorname{supp} f$ is the support of a function $f$
defined on $\Omega$, that is, $\operatorname{supp}
f:=\{\omega\in\Omega\dvtx f(\omega)\ne0\}$, $\mathcal{F}_\xi$ is
the distribution function of a random variable~$\xi$, and $[z]$
is the integral part of a real number $z$.

\subsection{The Kruglov property of symmetric function
spaces}\label{The Kruglov property}

Let $f$ be a
measurable function (a random variable) on $[0,1]$. By $\pi(f)$ we
denote the random variable $\sum_{i=1}^N f_i$, where $f_i$'s are
independent copies of $f$, and $N$ is a Poisson random variable with
parameter $1$ independent of the sequence $\{f_i\}$.
\begin{definition*}
A symmetric function space
$X$ is said to have the Kruglov property if and only if $f\in
X\Longleftrightarrow\pi(f)\in X$.
\end{definition*}

This property has been studied by Braverman
\cite{Br}, using some probabilistic constructions of
Kruglov \cite{K} and by the authors in \cite{AS1,AS2,AS3} via an
operator approach. We refer to the latter papers for various
equivalent characterizations of the Kruglov property. Note that
only the implication $f\in X\Longrightarrow\pi(f)\in X$ is
nontrivial, since the implication $\pi(f)\in X \Longrightarrow
f\in X$ is always satisfied \cite{Br}, page~11. Moreover, a symmetric
space $X$ has the Kruglov property if $X\supseteq L_p$ for some
$p<\infty$ \cite{Br}, Theorem 1.2, and \cite{AS1}, Corollaries 5.4,
5.6. At the same time, some exponential Orlicz spaces which
do not contain $L_q$ for any $q<\infty$ also possess this property
(see \cite{Br,AS1}).

\subsection{The Kruglov operator in symmetric function
spaces}\label{The Kruglov operator}

Let $\{B_n\}_{n=1}^\infty$ be a sequence of pairwise disjoint measurable
subsets of $[0,1]$ and let $\lambda(B_n)=\frac1{en!}$. If $f\in
L_1[0,1]$, then we set
%
%
\begin{equation}\label{Kruglov operator}
Kf(\omega_0,\omega_1,\ldots)=\sum_{n=1}^\infty\sum
_{k=1}^nf(\omega_k)\chi_{B_n}(\omega_0).
\end{equation}
Then $K\dvtx L_1[0,1]\to L_1(\Omega,{\mathcal P})$ is a positive linear
operator. Here $(\Omega,{\mathcal
P})=\prod_{n=0}^\infty([0,1],\lambda_n)$, where $\lambda_n$ is the
Lebesgue measure on $[0,1]$.

For convenience, by $Kf$ we also denote another random variable
defined on $[0,1]$ and having the same distribution as the variable
introduced in (\ref{Kruglov operator}). If $f\in L_1[0,1]$,
$\{B_n\}$ is the same sequence of subsets of $[0,1]$ as above, and,
for each $n\in\mathbb{N}$ $f_{n,1},f_{n,2},\ldots,f_{n,n}$ and $\chi_{B_n}$
form a set of independent functions such that $f_{n,k}^*=f^*$ for
every $k=1,\ldots,n$, then $Kf(t)$ is defined as the decreasing
rearrangement of the function
%
%
\begin{equation}\label{40-(-1)}
\sum_{n=1}^\infty\sum_{k=1}^n f_{n,k}(t)\chi_{B_n}(t)\qquad (0\le t\le1).
\end{equation}
As we have pointed out above $K$ is a linear operator from
$L_1[0,1]$ to $L_1(\Omega,{\mathcal P})$. By saying that $K$
maps boundedly a symmetric space $X$ on $[0,1]$ into symmetric
space~$Y$, we mean that $K$ is bounded as a linear mapping from
$X[0,1]$ into $Y(\Omega,{\mathcal P})$. The representation of
$Kf$ given by (\ref{40-(-1)}) allows us, without any ambiguity,
also speak about $K$ as a bounded map from $X[0,1]$ into
$Y[0,1]$. A direct computation (see, e.g., \cite{AS1}) yields the
following equality for the characteristic function $\varphi_{K
f}$ of $K f$:
%
%
\begin{eqnarray}\label{equalityfromAS1}
\varphi_{K f}(t)&=&\exp\biggl(\int_{-\infty}^\infty
(e^{itx}-1 )\,d\mathcal{F}_f(x) \biggr) \nonumber\\[-8pt]\\[-8pt]
&=& \exp\bigl(\varphi_f(t)-1\bigr)=\varphi_{\pi(f)}(t),\qquad
t\in{\mathbb R}.\nonumber
\end{eqnarray}
Therefore, $\mathcal{F}_{Kf}=\mathcal{F}_{\pi(f)}$, and we can treat
$Kf$ as an explicit representation of $\pi(f)$. In particular, a
symmetric space $X$ has the Kruglov property if and only if $K$ is
bounded in $X$.

It follows from the definition of the operator $K$ that for any
symmetric spaces $X$ and $Y$ $\|K\|_{X\to Y}\geq1/e$ provided
that $\|\chi_{[0,1]}\|_X=\|\chi_{[0,1]}\|_Y=1$ (see also
\cite{Br},  page 11). It is shown in \cite{AS1,AS2,AS3} that
the operator $K$ plays an important role in estimating the norm
of sums of i.r.v.'s through the norm of sums of their disjoint
copies. In particular, in \cite{AS1} the well-known results of
Johnson and Schechtman from \cite{JS} have been strengthened.
In the next section, we shall improve the main results of
\cite{AS1} and explain the role of this improvement in
obtaining sharp constants in Rosenthal-type inequalities
studied earlier in some special cases \cite{JS,JSZ,La}.
Subsequent sections contain explicit computation of the norm of
the operator $K$ in various classes of symmetric spaces $X$ and
further modifications of Rosenthal's inequality
(\ref{Rosenthal}).

\section{Kruglov operator and Rosenthal's inequalities}\label{The role}

The main objective of the present section is the strengthening of
\cite{AS1}, Theorem 3.5. For an arbitrary symmetric space $X$ on
$[0,1]$ and an arbitrary $p\in[1,\infty]$, we defined in
\cite{AS1} a function space $Z_X^p$ on $[0,\infty)$ by
\[
Z_X^p:=\{f\in L_1[0,\infty)+L_\infty[0,\infty)\dvtx \Vert
f\Vert'_{{Z_X^p}}<\infty\},
\]
where
\[
\Vert f\Vert'_{{Z_X^p}}:=\bigl\Vert f^*\chi_{{[0,1]}}\bigr\Vert
_{{X}}+\bigl\Vert
f^*\chi_{{[1,\infty)}}\bigr\Vert_{{p}}\asymp\bigl\Vert
f^*\chi_{{[0,1]}}\bigr\Vert_{{X}}+ \Biggl(\sum_{k=1}^\infty
f^*(k)^p \Biggr)^{1/p}.
\]
Clearly, $\Vert\cdot\Vert'_{{Z_X^p}}$ is a
quasi-norm. It is easy to see that $Z_X^p$ equipped with the
equivalent norm
\[
\Vert f\Vert_{{Z_X^p}}:=\bigl\Vert
f^*\chi_{{[0,1]}}\bigr\Vert_{{X}}+\Vert
f\Vert_{{(L_1+L_p)(0,\infty)}},\qquad f\in Z_X^p,
\]
is a symmetric space on $[0,\infty)$. The spaces $Z_X^2$ were
introduced in \cite{JMST}, and the spaces $Z_X^1$ and $Z_X^2$ were
used in \cite{JS} in the study of Rosenthal-type inequalities.
Following \cite{LT2}, page 46, we define the space $\widetilde
{X(l_p)}$ as the set of all sequences $\mathbf{f}=\{f_k(\cdot
)\}_{{k=1}}^\infty$, $f_k\in X$ $(k\ge1)$ such that
\[
\Vert\mathbf{f}\Vert_{{\widetilde{X(l_p)}}}:=\sup_{n=1,2,\ldots}
\Biggl\Vert\Biggl(\sum_{{k=1}}^n|f_k|^p \Biggr)^{1/p} \Biggr\Vert
_{{X}}<\infty
\]
(with an obvious modification for $p=\infty$). The closed subspace of
$\widetilde{X(l_p)}$ generated by all eventually vanishing sequences
$\mathbf{f}\in\widetilde{X(l_p)}$ is denoted by $X(l_p)$.

Let $X$ and $Y$ be symmetric spaces on $[0,1]$ such that
$X\subseteq Y$. The main focus of \cite{AS1,AS2,AS3} is on
inequalities of the type
%
%
\begin{equation}\label{maininequality} \Biggl\Vert
\sum_{i=1}^nf_i \Biggr\Vert_{{Y}}\leq C \Biggl\Vert\sum
_{i=1}^n\bar f_i \Biggr\Vert_{{Z_X}}
\end{equation}
and
%
%
\begin{equation}\label{maininequality'}\Vert\mathbf{f}\Vert
_{{Y(l_p)}}\leq
C \Biggl\Vert\sum_{i=1}^n\bar f_i \Biggr\Vert_{{Z_X^p}},
\end{equation}
where
the sequence $\mathbf{f}:=\{f_k\}_{k=1}^\infty\subset X$ consists
of i.r.v.'s, and the sequence $\{\bar f_k\}_{k=1}^\infty$ is a
disjointly supported sequence of equimeasurable copies of the
elements from the sequence $\mathbf{f}$.

Our first result in this section strengthens \cite{AS1},
Theorems 3.5 and 6.1, by establishing sharp estimates on the constant
$C$ in
(\ref{maininequality}). Let $F$ be the disjointification function
related to the sequence $\mathbf{f}:=\{f_k\}_{k=1}^\infty$ [see
(\ref{f1})].
\begin{theorem}\label{precise1} Let $X$ and $Y$ be symmetric spaces on
$[0,1]$ such that
$X\subseteq Y$ and $Y$ has the Fatou property.

\begin{longlist}
\item If there exists a constant $C$ such that the
estimate
%
%
\begin{equation}\label{maininequality1} \Biggl\Vert\sum
_{i=1}^nf_i \Biggr\Vert_{{Y}}\leq C\Vert F\Vert
_{Z_X}\qquad \Biggl(=C \Biggl\Vert\sum_{i=1}^n\bar f_i \Biggr\Vert
_{{X}} \Biggr),
\end{equation}
holds for every sequence
$\{f_k\}_{k=1}^\infty\subset X$ of i.r.v.'s satisfying
the assumption (\ref{assumption}) for all $n\in\mathbb{N}$,
then the operator $K$ acts boundedly from
$X$ into $Y$ and $\|K\|_{X\to Y}\leq C$.
\item If the operator $K$ acts boundedly from
$X$ into $Y$, then for every sequence $\{f_k\}_{k=1}^\infty\subset
X$ of independent random variables, we have
%
%
\begin{equation}\label{maininequality2} \Biggl\Vert\sum
_{i=1}^nf_i \Biggr\Vert
_{{Y}}\leq\alpha\|K\|_{X\to Y}\Vert F\Vert_{{Z_X}},
\end{equation}
where
$\alpha>0$ is a universal constant which does not depend on $X$ and
$Y$.
\end{longlist}
\end{theorem}
\begin{pf}
(i) The claim follows from the inspection of the first part of the
proof of \cite{AS1}, Theorem 3.5. For the convenience of the reader, we
include details of the argument. Fix $f\in X$ and $n\in\mathbb{N}$ and
choose $h\in X$ such that $\mathcal{F}_h=\mathcal{F}_f$ and such
that $h$ and $\chi_{{[0,1/n]}}$ are independent. Set
$h_n:=h\chi_{{[0,1/n]}}$, and let $\{\chi_{{[0,1/n]}},
h_{{n,k}}\}_{{k=1}}^n$ be a set of $(n+1)$ independent random
variables such that $\mathcal{F}_{h_{{n,k}}}=\mathcal{F}_{h_n}$
for all $1\leq k\leq n$. Since the functions $ |{\sum_{{k=1}}^n
\bar h_{{n,k}}} |$ and $|h|$ have the same distribution
function, we conclude that the functions $ |{\sum_{{k=1}}^n
\bar h_{{n,k}}} |$ and $|f|$ are equidistributed. Observing
now that the sequence $\{ h_{{n,k}}\}_{k=1}^n$ satisfies
(\ref{assumption}), we obtain
%
%
\begin{equation}\label{3.2fromAS1}
\Biggl\Vert\sum_{{k=1}}^n h_{{n,k}} \Biggr\Vert_{{Y}}\leq
C \Biggl\Vert\sum_{{k=1}}^n \bar h_{{n,k}} \Biggr\Vert_{{X}}=C\Vert
f\Vert_{{X}}.
\end{equation}

A direct computation shows that
$\varphi_{{h_n}}(t)=n^{-1}\varphi_f(t)+(1-n^{-1})$ for all $t\in
{\mathbb R}$. Hence, the characteristic function of the sum
$H_n:=\sum_{{k=1}}^n h_{{n,k}}$ is given by
\[
\varphi_{{H_n}}(t)= \bigl(n^{-1}\bigl(\varphi_f(t)-1\bigr)+1
\bigr)^n\qquad
\forall t\in{\mathbb R}.
\]
Since $\lim_{{n\to
\infty}}\varphi_{{H_n}}(t)=\exp(\varphi_f(t)-1)=\varphi_{{\pi(f)}}(t)$,
for all $t\in{\mathbb R}$, we see that $H_n$ converges weakly to $K
f$. Combining this with (\ref{3.2fromAS1}), (\ref{equalityfromAS1}),
\cite{Br}, Proposition 1.5, and with the fact that $Y$ has the Fatou
property, we conclude that
\[
\Vert K f\Vert_{{Y}}\leq C\Vert
f\Vert_{{X}}.
\]
This completes the proof of the first assertion.

(ii) Firstly, let us assume that a sequence
$\{f_k\}_{k=1}^\infty\subset X$ consists of independent random
variables satisfying assumption (\ref{assumption}). Denote by
$\{h_k\}_{{k=1}}^n$ a sequence of independent random variables
such that $\mathcal{F}_{h_k}=\mathcal{F}_{{\pi(f_k)}}$ for all
$k=1,2,\ldots.$ Consider the following two cases.

(a) $f_k$'s are symmetrically distributed r.v.'s. In \cite{Pr},
Prokhorov proved that in this case we have
\[
\lambda\Biggl\{ \Biggl|\sum_{{k=1}}^nf_k \Biggr|\ge x \Biggr\}\leq
8\lambda\Biggl\{ \Biggl|\sum_{{k=1}}^nh_k \Biggr|\ge\frac{x}{
2} \Biggr\} \qquad(x>0).
\]
From this inequality (see, e.g., \cite{KPS}, Corollary II.4.2) it
follows that
%
%
\begin{equation}\label{fromProh}
\Biggl\|\sum_{{k=1}}^n f_k \Biggr\|_{{Y}}\leq16 \Biggl\|\sum_{{k=1}}^n
h_k \Biggr\|_{{Y}}.
\end{equation}

(b) $f_k=a_k\chi_{A_k}$, where $a_k\ge0$ and $A_k$ are arbitrary
independent subsets of $[0,1]$ $(k=1,2,\ldots)$. Without loss of
generality, we may assume that $f_k$'s are defined on the measure
space $\prod_{n=0}^\infty([0,1],\lambda_n)$ by the formula
$f_k(t)=a_k\chi_{[0,p_k]}(t_k)$, where $p_k=\lambda(A_k)$, $k\ge
1$. From the definition of the Kruglov operator [see
(\ref{40-(-1)})] it follows that $\lambda(\{Kf_k=a_k\})=p_k/e\geq
p_k/3$. Hence, by (\ref{equalityfromAS1}), we may assume that
$h_k\geq a_k\chi_{[0,p_k/3]}$, and so
%
%
\begin{equation}\label{fromDEFIN}
\Biggl\|\sum_{k=1}^nf_k\Biggr\|_Y\leq
3\Biggl\|\sum_{k=1}^na_k\chi_{[0,p_k/3]}(t_k)\Biggr\|_Y\leq3
\Biggl\|\sum_{k=1}^nh_k\Biggr\|_Y.
\end{equation}

Next, from $\bar f_k\bar f_m=0$ $(k\neq m)$ it follows that
\[
e^{itF}-1=\sum_{k=1}^n(e^{it\bar f_k}-1).
\]
Therefore,
\[
\varphi_F(t)-1=\int(e^{itF}-1)=\sum_{k=1}^n\int(e^{itf_k}-1)=\sum
_{k=1}^n\bigl(\varphi_{f_k}(t)-1\bigr)
\]
and
\[
\varphi_{K(F)}=\exp(\varphi_F-1)=\prod_{k=1}^n\exp(\varphi
_{f_k}-1)=\prod_{k=1}^n\varphi_{h_k}
=\varphi_{\sum_{k=1}^nh_k}.
\]
Thus, the sum $\sum_{k=1}^nh_k$ is equidistributed with $K(F)$.
Therefore, by inequalities (\ref{fromProh}) and (\ref{fromDEFIN}),
we have
%
%
\begin{equation}\label{relation 1}
\Biggl\|\sum_{k=1}^n f_k \Biggr\|_Y \leq16 \Biggl\|\sum_{k=1}^n
h_k \Biggr\|_Y=16\|K(F)\|_Y\leq16\|K\|_{X\to Y}\|F\|_X
\end{equation}
in the case (a) and analogously
%
%
\begin{equation}\label{relation 2}
\Biggl\|\sum_{k=1}^n f_k \Biggr\|_Y \leq3\|K\|_{X\to Y}\|F\|_X
\end{equation}
in the case (b).

Now, we consider the case when a sequence
$\{f_k\}_{k=1}^\infty\subset X$ consists of mean zero i.r.v.'s
satisfying assumption (\ref{assumption}) with $1/2$ instead of
$1$. We shall use the standard ``symmetrization trick.'' Let
$\{f_k'\}_{k=1}^n$ be a sequence of i.r.v.'s such that the
sequence $\{f_k,f_k'\}_{k=1}^n$ consists of independent random
variables and $\mathcal{F}_{f_k}=\mathcal{F}_{f_k'}$ for all $k\ge
1$.

Setting $g_k:=f_k-f_k'$, we obtain a sequence $\{g_k\}_{{k=1}}^n$
of symmetrically distributed independent random variables
satisfying assumption (\ref{assumption}), and, by (\ref{relation 1}),
we have
%
%
\begin{equation}\label{maininequality2H} \Biggl\Vert\sum_{k=1}^n
g_k \Biggr\Vert_Y\leq16 \|K\|_{X\to Y}\Vert G\Vert_X,
\end{equation}
where
$G:=\sum_{{k=1}}^n\bar g_k$. Let $\mathcal{B}$ be the
$\sigma$-subalgebra generated by the sequence
$\{f_k\}_{{k=1}}^n$ and let $E_{\mathcal{B}}$ be the
corresponding conditional expectation operator. Thanks to our
assumption $Y$ is an 1-interpolation space for the couple
$(L_1,L_\infty)$. Hence, since $E_{\mathcal{B}}$ is bounded in
$L_1$ and $L_\infty$ (with constant $1$)
\cite{LT2}, Theorem 2.a.4, we have $\|E_{\mathcal{B}}\|_{Y\to
Y}=1$. Therefore, due to the independence of $f_k$'s and $f_k'$'s,
we have
%
%
\begin{equation}\label{maininequality2''} \Biggl\Vert\sum_{k=1}^n
f_k \Biggr\Vert_Y= \Biggl\Vert E_{\mathcal{B}} \Biggl(\sum_{k=1}^n
g_k \Biggr) \Biggr\Vert_Y\leq\Biggl\Vert\sum_{k=1}^n g_k \Biggr\Vert_Y.
\end{equation}
On the other hand, it is obvious that $\|G\|_X\leq2\|F\|_X$.
Combining (\ref{maininequality2''}) and (\ref{maininequality2H}),
we see that (\ref{maininequality2}) holds with $\alpha=32$.

Next, let us consider an arbitrary sequence
$\{f_k\}_{k=1}^\infty\subset X$ of i.r.v.'s satisfying assumption
(\ref{assumption}) with $1/2$ instead of $1$. In this case, set
$u_k=f_k-v_k$, where
\[
v_k:=\frac1{\lambda(\{\operatorname{supp} f_k\})}\int_0^1f_k(t) \,dt\cdot\chi
_{\operatorname{supp} f_k}\qquad (k\ge1).
\]
Clearly, $\{u_k\}_{k\geq1}$ is a mean zero sequence of i.r.v.'s
satisfying assumption (\ref{assumption}) with $1/2$ instead of
$1$. Thus, the preceding argument yields
%
%
\begin{equation}\label{maininequality21} \Biggl\Vert\sum_{k=1}^n
u_k \Biggr\Vert_Y\leq32 \|K\|_{X\to Y}\Vert U\Vert_X,
\end{equation}
where
$U:=\sum_{k=1}^n\bar u_k$. Moreover, by (\ref{relation 2}), we
have
\[
\Biggl\Vert\sum_{k=1}^n f_k \Biggr\Vert_Y\leq\Biggl\Vert\sum_{k=1}^n u_k
\Biggr\Vert_Y+ \Biggl\Vert\sum_{k=1}^n v_k \Biggr\Vert_Y\leq32 \|K\|
_{X\to
Y}\Vert U\Vert_X+3\|K\|_{X\to Y}\Vert V\Vert_X.
\]
Let $\mathcal{C}$ be the $\sigma$-algebra generated by the
supports of $\bar f_k$'s. It is clear that $V=E_{\mathcal{C}}(F)$,
where $V:=\sum_{k=1}^n\bar v_k$. Therefore, as above, we have
$\|V\|_X\leq\|F\|_X$. Moreover, since $U=F-V$, we also have
$\|U\|_X\leq\|F\|_X+\|V\|_X\leq2\|F\|_X$. Thus,
\[
\Biggl\Vert\sum_{k=1}^n f_k \Biggr\Vert_Y\leq67\|K\|_{X\to Y}\Vert
F\Vert_X.
\]

If $\{f_k\}_{k=1}^\infty\subset X$ is an arbitrary sequence of
i.r.v.'s satisfying assumption (\ref{assumption}), then we may
represent $f_k=f_k'+f_k''$, where each of the sequences
$\{f_k'\}_{k=1}^\infty$ and $\{f_k''\}_{k=1}^\infty$ consists of
i.r.v.'s satisfying assumption (\ref{assumption}) with $1/2$
instead of $1$ and moreover $|F'|\le|F|$ and $|F''|\le|F|$. In
this case, using the preceding formulas, we obtain that
\[
\Biggl\Vert\sum_{k=1}^n f_k \Biggr\Vert_Y\leq\Biggl\Vert\sum_{k=1}^n f_k'
\Biggr\Vert_Y+ \Biggl\Vert\sum_{k=1}^n f_k'' \Biggr\Vert_Y\leq134
\|K\|_{X\to Y}\Vert F\Vert_X.
\]

Finally, repeating verbatim the proof of \cite{AS1}, Theorem 6.1,
we obtain (\ref{maininequality2}) for arbitrary sequences of
i.r.v.'s [which do not necessarily satisfy assumption
(\ref{assumption})] as a corollary of already considered special
case when (\ref{assumption}) holds.
\end{pf}

Our next result strengthens \cite{AS1}, Theorem 6.7, by establishing
sharp estimates on the constant $C$ in (\ref{maininequality'}).
Before proceeding, we recall the following construction due to
Calderon \cite{Cal}. Let $X_0$ and $X_1$ be two Banach lattices of
measurable functions on the same measure space $(\mathcal{M},m)$
and let $\theta\in(0,1)$. The space $X_0^{1-\theta}X_1^{\theta}$
consists of all measurable functions $f$ on $(\mathcal{M},m)$ such
that for some $\lambda
>0$ and $f_i\in X_i$ with $\Vert f_i\Vert_{{X_i}}\leq1$,
$i=0,1$, we have
\[
|f(x)|\leq\lambda|f_0(x)|^{1-\theta}|f_1(x)|^{\theta},\qquad x\in
\mathcal{M}.
\]
This space is equipped with the norm given by the greatest lower
bound of all numbers $\lambda$ taken over all possible such
representations. Even though this construction is not an
interpolation functor on general couples of Banach lattices
(see~\cite{Lo}), it is still a convenient tool of interpolation theory.
Indeed, if $(X_0,X_1)$ is a Banach couple and if ($Y_0,Y_1)$ is
another Banach couple of lattices of measurable functions on some
measure space $(\mathcal{M}',m')$, then any positive operator $A$
from $S(\mathcal{M},m)$ into $S(\mathcal{M}',m')$, which acts
boundedly from the couple ($X_0,X_1)$ into the couple ($Y_0,Y_1)$
also maps boundedly $X_0^{1-\theta}X_1^{\theta}$ into
$Y_0^{1-\theta}Y_1^{\theta}$ and, in addition, $\Vert A\Vert
_{{X_0^{1-\theta}X_1^{\theta}\to Y_0^{1-\theta}Y_1^{\theta}}}\leq
\Vert A\Vert^{1-\theta} _{{X_0\to Y_0}}\Vert A\Vert^{\theta}
_{{X_1\to Y_1}}$ for all $\theta\in(0,1)$. The proof of the latter
claim follows by inspection of the standard arguments from
\cite{LT2}, Proposition~1.d.2(i), page 43.
\begin{theorem}\label{precise2} Let $X$ and $Y$ be symmetric spaces on
$[0,1]$ such that
$X\subseteq Y$. If $K$ acts boundedly from $X$ into $Y$, then there
exists a universal constant $\beta>0$ such that for every sequence
$\{f_k\}_{k=1}^\infty\subset X$ of i.r.v.'s and for every $q\in
[1,\infty]$, we have
\[
\Vert\mathbf{f}\Vert_{Y(l_q)}\leq\beta\|K\|_{X\to Y}^{1/q}\Vert F\Vert
_{Z_X^q}.
\]
\end{theorem}
\begin{pf} The case $q=1$ has been treated in Theorem \ref{precise1}. If
$q=\infty$, then it is sufficient to observe that
\[
\Vert F^*\Vert_{Z_X^\infty}\asymp\bigl\Vert F^*\chi_{[0,1]}\bigr\Vert_{X},
\]
and that (see, e.g., \cite{HM}, Proposition 2.1)
%
%
\begin{equation}\label{Proposition 2.1HM} \frac{1}{2}\lambda\bigl\{
F^*\chi_{{[0,1]}}>\tau\bigr\}\leq\lambda\Bigl\{{\sup_{{k=1,2,\ldots
}}}|f_k|>\tau\Bigr\}\leq\lambda\bigl\{
F^*\chi_{{[0,1]}}>\tau\bigr\}\qquad
(\tau>0).\hspace*{-25pt}
\end{equation}
Therefore, for some constant $\gamma$ (which does not
depend on $X$ and $Y$) we have
\[
\Vert\mathbf{f}\Vert_{{Y(l_\infty)}}\leq\Vert\mathbf{f}\Vert_{
{X(l_\infty)}}\leq\gamma\Vert F\Vert_{{Z_X^\infty}}.
\]
The rest of the proof is based on methods from interpolation theory
and is very similar to the arguments in \cite{AS1}, Theorem 6.7. Let
$\delta\dvtx (\Omega, {\mathcal P})\to([0,1],\lambda)$ be a measure
preserving isomorphism, where $(\Omega, {\mathcal
P}):=\prod_{{k=0}}^\infty([0,1], \lambda_k)$ (here, $\lambda_k$ is
the Lebesgue measure on $[0,1]$ for every $k\ge0$). For every $g\in
S(\Omega, {\mathcal P})$, we set $T(g)(x):=g(\delta^{-1}x)$, $x\in
[0,1]$. Note that $T$ is a rearrangement-preserving mapping between
$S(\Omega, {\mathcal P})$ and $S([0,1],\lambda)$. We define the
positive linear mapping $Q$ from $S(0,\infty)$ into $S(\Omega,
{\mathcal P})^{\mathbb{N}\cup\{0\}}$ by setting
\[
Qf(\omega_0,\omega_1,\ldots):=\{(Qf)_k\}_{{k=0}}^\infty,\qquad f\in
S(0,\infty),
\]
where $(Qf)_k(\omega_0,\omega_1,\ldots):=f(\omega_k+k)$ $(\omega
_k\in
[0,1])$, $k\ge0$. The arguments above show that the positive
operator $Q'f:=\{T(Qf)_k\}_{{k=0}}^\infty$ is bounded from $Z_X^1$
into $Y(l_1)$ with the norm not exceeding $\alpha\|K\|_{X\to Y}$
(where $\alpha$ is a universal constant) and also from $Z_X^\infty$
into $Y(l_\infty)$ with the norm not exceeding $\gamma$. It follows
that
\[
Q'\dvtx (Z_X^1)^{1-\theta}(Z_X^\infty)^{\theta}\to
(Y(l_1))^{1-\theta}(Y(l_\infty))^{\theta}
\]
and
\[
\|Q'\|\leq\|Q'\|_{Z_X^1\to Y(l_1)}^{1-\theta}\cdot\|Q'\|_{Z_X^\infty
\to Y(l_\infty)}^{\theta},\qquad \theta\in(0,1).
\]
Now, we shall use the following facts: for every $\theta\in(0,1)$
we have
\[
Z_X^q\subseteq(Z_X^1)^{1-\theta}(Z_X^\infty)^{\theta},\qquad
(Y(l_1))^{1-\theta}(Y(l_\infty))^{\theta}\subseteq Y(l_q),\qquad
q=\frac{1}{ {1-\theta}}.
\]
To see the first embedding above, fix $g=g^*\in Z_X^q$, $\Vert
g\Vert_{{Z_X^q}}=1$, and set
\[
g_1:=g\chi_{{[0,1]}}+g^q\chi_{{[1,\infty)}},\qquad
g_\infty:=g\chi_{{[0,1]}}+\chi_{{[1,\infty)}}.
\]
Clearly, $g=(g_1)^{1-\theta}(g_\infty)^{\theta}$. Moreover, since
$g(1)\le1$, then $g_1$ decreases, which implies that $g_i\in Z_X^i$
and $\Vert g_i\Vert_{{ Z_X^i}}\leq3$ $(i=1,\infty)$. The second
embedding above (in fact, equality) is established in \cite{Bu},
Theorem 3. Now, we are in a position to conclude that
\[
\|Q'\|_{Z_X^q\to Y(l_q)}\leq\alpha^{1/q}\|K\|_{X\to Y}^{1/q}\cdot
\gamma^{1-1/q}\leq\max(1,\alpha,\gamma)\|K\|_{X\to Y}^{1/q}.
\]
The proof is completed by noting that every sequence
$\{f_k\}_{k=1}^\infty\subset X$ of i.r.v.'s may be represented in
the form $\{f_k\}= Q'(\tilde{F})$, with some function $\tilde{F}$
which is equidistributed with $F$.
\end{pf}

The results presented in this section show that the sharp estimates in
the deterministic estimates of expressions
\[
\Biggl\| \sum
_{{k=1}}^nf_k \Biggr\|_X, \Biggl\| \Biggl(\sum
_{{k=1}}^n|f_k|^p \Biggr)^{1/p} \Biggr\|_X,\qquad 1\leq p<\infty,
\]
in
terms of the sum of disjoint copies of individual terms of
$\mathbf{f}$ are fully determined by the norm of the Kruglov
operator $K$ in $X$. In the next section, we shall present
sharp estimates of this norm in a number of important cases,
including the case $X=L_p$ $(1\leq p<\infty)$ studied earlier
in \cite{JSZ,La} and \cite{Ju} by completely different
methods. It does not seem that the methods used in those papers
can be extended outside of the $L_p$-scale.

\section{Norm of the Kruglov operator and sharp constants in
Rosenthal's inequalities}\label{The norm}

Recall that a Banach lattice $X$ is said to satisfy an
\textit{upper} $p$-\textit{estimate}, if there exists a constant
$C>0$ such that for every finite sequence
$(x_{j})_{j=1}^{n}\subseteq X$ of pairwise disjoint elements,
\[
\Biggl\Vert\sum_{j=1}^{n}x_{j} \Biggr\Vert_{X}\leq C \Biggl(
\sum_{j=1}^{n}\Vert x_{j}\Vert_{X}^{p} \Biggr) ^{1/p}.
\]

Recall also that if $\tau>0$, \textit{the dilation operator}
$\sigma_{\tau}$ is defined by setting
\[
(\sigma_{\tau}x)(s)=
\cases{
x(s/\tau),&\quad $s\leq\min\{1,\tau\}$,\cr
0,&\quad $\tau<s\leq1$.}
\]
The operator $\sigma_{\tau}$, $\tau>0$ acts boundedly in every
symmetric function space $X$ \cite{KPS}, Theorem II.4.4.

First, we suppose that $1\leq p<\infty$ and that $X$ is a symmetric
space satisfying the following two conditions:
\begin{longlist}
\item $X$ satisfies an upper $p$-estimate;
\item $\|\sigma_t\|_{X\to X}\leq Ct ^{1/p},
0<t<1$.
\end{longlist}
Conditions (i) and (ii) imply, in particular, that both Boyd
indices of $X$ (see, e.g., \cite{LT2}) are equal to $p$.

The value of the constant $C$ varies from line to line in this
section.
\begin{proposition}\label{estimateLp} If a symmetric space $X$
satisfies the above assumptions,
then there exists a universal constant $\alpha>0$ whose value
depends only on the constants in \textup{(i)} and \textup{(ii)}
above such that
%
%
\begin{equation}\label{mainLp}
\|K\|_{X\to X}\leq\alpha\frac{p}{
{\ln(p+1)}},\qquad
p\ge1.
\end{equation}
\end{proposition}
\begin{pf}
Let $0\leq f\in X$, $p\ge1$, and $n\in\mathbb{N}$. Let $f_{n, 1},
f_{n,2},\ldots, f_{n,n}$ and $\chi_{B_n}$ have the same meaning
as in Section \ref{The Kruglov operator}.
Then, by (\ref{40-(-1)}), the random variable $Kf$ is equimeasurable
with the random variable
\[
\sum_{n=1}^\infty g_{n}\chi_{B_n}\qquad \mbox{where } g_n=\sum
_{k=1}^n f_{n,k} \qquad(n=1,2,\ldots).
\]
Since $g_n$ and $\chi_{B_n}$ are independent, then the
assumptions on $X$ imply
\begin{eqnarray*}
\|Kf\|_X & \leq& C \Biggl(\sum
_{n=1}^\infty\|g_n\chi_{B_n}\|_X^p \Biggr)^{1/p}= C \Biggl(\sum
_{n=1}^\infty\bigl\|\sigma_{\lambda(B_n)}g_n\bigr\|_X^p \Biggr)^{1/p}\\
&\leq&
C \Biggl(\sum_{n=1}^\infty\frac{1}{e\cdot
n!}\|g_n\|_X^p \Biggr)^{1/p}\le C \Biggl(\sum_{n=1}^\infty
\frac{n^p}{n!} \Biggr)^{1/p}\|f\|_X.
\end{eqnarray*}

It is clear that
\[
\sum_{n=1}^{\infty}\frac{n^p}{n!}\leq\sup_n\frac{n^p}{p^n}\sum
_{n=1}^{\infty}\frac{p^n}{n!}=e^p\sup_{n}\frac{n^p}{p^n}.
\]
The function $x^p/p^x$ takes its maximal value at $x=p/\log(p)$ and
this maximum does not
exceed $(p/\log(p))^p$. Therefore,
\[
\Biggl(\sum_{n=1}^{\infty}\frac{n^p}{n!}\Biggr)^{1/p}\leq e\cdot\frac{p}{\log(p)}.
\]
On the other side, if $1\le p<2$, then
\[
\Biggl(\sum_{n=1}^{\infty}\frac{n^p}{n!} \Biggr)^{1/p}\leq\sum
_{n=1}^{\infty}\frac{n^2}{n!}=2e.
\]
\upqed\end{pf}

It is well known that an $L_p$-space, $1\leq p<\infty$, satisfies
an upper $p$-estimate and that $\|\sigma_t\|_{L_p\to L_p}=t
^{1/p}, 0<t\le1$. The facts that similarly $M(t^{1/p})$,
$1<p<\infty$, satisfies an upper $p$-estimate and
$\|\sigma_t\|_{M(t^{1/p})\to M(t^{1/p})}\leq t ^{1/p},
0<t<1$, follow from a combination of \cite{Cr}, Theorem 3.4(a)(i), and
\cite{LT2}, Proposition 1.f.5, and from~\cite{KPS}, Chapter II, Theorem 4.4.
Combining these facts with Proposition \ref{estimateLp}, we obtain
the following corollary.
\begin{corollary}\label{corollaryLp} There exists a constant $C>0$
such that for all $p\ge1$
%
%
\begin{equation}\label{mainLpcorollary} \|K\|_{L_p\to L_p}\leq\frac{Cp}{
{\ln(p+1)}} \quad\mbox{and}\quad \|K\|_{M(t^{1/p})\to
M(t^{1/p})}\leq\frac{Cp}{ {\ln(p+1)}}.
\end{equation}
\end{corollary}

Although the Lorentz space $\Lambda(t^{1/p})$, $1<p<\infty$,
does not satisfy the assumptions of Proposition
\ref{estimateLp}, nevertheless estimates similar to
(\ref{mainLpcorollary}) also hold for the norm of the operator
$K\dvtx\Lambda(t^{1/p})\to\Lambda(t^{1/p})$, $1\leq p<\infty$. The
proof below is based on the properties of the Kruglov operator
$K$ in Lorentz spaces exposed in \cite{AS1}, Section~5.
\begin{proposition}\label{mainLambdapproposition} There exists a
constant $C>0$ such that for all $p\ge1$
%
%
\begin{equation}\label{mainLambdap}
\|K\|_{\Lambda(t^{1/p})\to\Lambda(t^{1/p})}\leq\frac{Cp}{ {\ln
(p+1)}},\qquad p\ge1.
\end{equation}
\end{proposition}
\begin{pf} By \cite{AS1}, Theorem 5.1, we have
\[
\|K\|_{\Lambda(\psi)\to\Lambda(\psi)}\leq2\sup_{u\in(0,1)}\frac
1{\psi(u)}\sum_{k=1}^{\infty}\psi\biggl(\frac{u^k}{k!}\biggr).
\]
If $\psi(t)=t^{1/p}$, then the latter supremum is equal to
\[
\sum_{n=1}^{\infty} \biggl(\frac1{n!} \biggr)^{1/p}\leq2+\sum
_{n=3}^{\infty}
\biggl(\frac1{[{n}/3]} \biggr)^{3[{n}/{3}]/p}=2+3\sum
_{k=1}^{\infty}e^{-3k\log(k)/p}.
\]
Since the latter sequence decreases, we can replace sum with an
integral and obtain
\[
\|K\|_{\Lambda_{t^{1/p}}\to\Lambda_{t^{1/p}}}\leq22+6\int
_e^{\infty}e^{-3s\log(s)/p}\,ds.
\]
Substitute $t=s\log(s)$. It follows that
$\frac{dt}{ds}=1+\log(s)\geq\frac12\log(t)$. Therefore,
\[
\|K\|_{\Lambda_{t^{1/p}}\to\Lambda_{t^{1/p}}}\leq22+12\int
_e^{\infty}\frac1{\log(t)}e^{-3t/p}\,dt.
\]
It is clear that
\[
\int_e^p\frac1{\log(t)}e^{-3t/p}\,dt\leq\int_e^p\frac1{\log
(t)}\,dt\leq \operatorname{const}\cdot\frac{p}{\log(p)}.
\]
On the other side,
\[
\int_p^{\infty}\frac1{\log(t)}e^{-3t/p}\,dt\leq\frac1{\log(p)}\int
_p^{\infty}e^{-3t/p}\,dt\leq\frac{p}{3\log(p)}.
\]
It follows that
\[
\|K\|_{\Lambda_{t^{1/p}}\to\Lambda_{t^{1/p}}}\leq \operatorname{const}\cdot\frac
{p}{\log(p)}.
\]
\upqed\end{pf}

We shall now estimate the norm of the operator
$K\dvtx\Lambda(t^{1/p})\to M(t^{1/p})$, $1< p<\infty$.
\begin{lemma}\label{mainweakestimate}
%
%
\begin{equation}\label{mainLambdaMp}
\|K\|_{\Lambda(t^{1/p})\to M(t^{1/p})}\asymp\frac{p}{\ln(p+1)},\qquad
1< p<\infty.
\end{equation}
\end{lemma}
\begin{pf} Since $\Lambda(t^{1/p})\subset L_p\subset M(t^{1/p})$ for
all $1<p<\infty$, the estimate from
the above in (\ref{mainLambdaMp}) follows immediately from Corollary
\ref{corollaryLp}. Let us now concentrate on the converse
inequality. Since $\|K1\|_{M(t^{1/p})}\ge1-1/e$ for all $p\ge1$,
then it suffices to establish the estimate from below for
sufficiently large $p$'s.

By \cite{AS1}, Remark 5.2, we have
%
%
\begin{eqnarray}\label{mainLambdaMp1}
\|K\|_{\Lambda(t^{1/p})\to M(t^{1/p})} &\ge&
\frac1e \sup_{u\in(0,1], k\in\mathbb{N}}\frac{k (
{u^k}/{k!} )^{1/p}} {u^{1/p}}=\frac1e \sup_{k\in\mathbb{N}}\frac
{k}{(k!)^{1/p}}\nonumber\\[-8pt]\\[-8pt]
&\geq&\frac1e\sup_{x\in\mathbb{R}}\frac{x}{x^{x/p}}.\nonumber
\end{eqnarray}
Substitute $x=p/\log(p)$. We obtain
\[
\|K\|_{\Lambda(t^{1/p})\to M(t^{1/p})}\ge\frac1e\bigl(p/\log
(p)\bigr)^{1-1/\log(p)}\geq\frac1{e}\frac{p}{\log(p)}.
\]
\upqed\end{pf}

The first main result of the present section is the following theorem
and its corollary.
\begin{theorem}\label{Lpnorm} We have
%
%
\begin{eqnarray}\label{Lpnormasymp}
\|K\|_{L_p\to L_p}&\asymp&\|K\|_{\Lambda(t^{1/p})\to\Lambda
(t^{1/p})}\asymp\|K\|_{M(t^{1/p})\to M(t^{1/p})}\nonumber\\[-8pt]\\[-8pt]
&\asymp&\frac{p}{\ln
(p+1)},\qquad p> 1,\nonumber
\end{eqnarray}
with universal constants.
\end{theorem}
\begin{pf}
The proof follows from combining Proposition \ref
{mainLambdapproposition} and Corollary \ref{corollaryLp} with Lemma
\ref{mainweakestimate}.
\end{pf}
\begin{corollary}\label{Rosenthal1}
The order of the constant $\frac{\alpha p}{\ln(p+1)}$ in
Rosenthal's inequality (\ref{Rosenthal}) with $X=L_p$ is optimal
when $p\to\infty$.
\end{corollary}
\begin{pf}
Apply Theorem \ref{precise1} to the case $X=Y=L_p$ and then apply
Theorem~\ref{Lpnorm}.
\end{pf}
\begin{remark} The same argument as above also shows that the order of
the constant
$\frac{\alpha p}{\ln(p+1)}$ is optimal in variants of Rosenthal's
inequality (\ref{Rosenthal}) for scales of Lorentz spaces
$\Lambda(t^{1/p})$ ($1< p<\infty)$ and Marcinkiewicz spaces
$M(t^{1/p})$ ($1< p<\infty)$.
\end{remark}
\begin{remark}
Earlier the result presented in Corollary \ref{Rosenthal1} was
established in \cite{JSZ,La} and \cite{Ju} by completely
different methods. Our approach here shows that the order of the
constant whether in the special $L_p$-case studied in the papers
just cited, or in a more general case of various scales of
symmetric spaces (as indicated in the preceding remark) is
fully determined by the norm of the Kruglov operator.
\end{remark}
\begin{remark}\label{meanzerocase}
It is shown in \cite{JSZ}, Theorem 4.1 and Proposition 4.3, that
the constant $\frac{\alpha p}{1+\ln(p)}$, where $\alpha$ is an
absolute constant, is also sharp in order when $p\to\infty$ in
Rosenthal's inequality (\ref{meanRosenthal}) for mean zero
i.r.v.'s. More precisely, using our notation, it is proved
there that for any such sequence ${\{f_k\}_{k=1}^n}\subset L_p$ we have
%
%
\begin{equation}\label{LpRosenthal}
\Biggl\|\sum_{k=1}^n f_k\Biggr\|_p\leq\frac{\alpha p}{1+\ln(p)}
\Biggl(\bigl\|F^*\chi_{[0,1]}\bigr\|_p+ \Biggl(\sum_{k=1}^nF^*(k)^2 \Biggr)^{1/2} \Biggr).
\end{equation}
Furthermore, if $c_p$ is the least constant in similar
inequality which would hold for any sequence ${\{f_k\}_{k=1}^n}\subset L_p$
of symmetrically and identically distributed i.r.v.'s, then
$c_p\ge\frac{p}{\sqrt{2} e(1+\ln(p))}$. Without going into
precise details, we observe that a careful inspection of the
proof of \cite{AS3}, Theorem 3.1, shows that a similar result
holds also in the case of an arbitrary symmetric space $X$ with
Kruglov's property (with an obvious replacement of the constant
$\frac{\alpha p}{1+\ln(p)}$ with the constant $\|K\|_{X\to
X}$). We state this result in full.
\begin{proposition}\label{meanzeroXcase} There exists an absolute
constant $\alpha$ such that
if the Kruglov operator $K$ acts boundedly in a symmetric space $X$,
then for any $n\in N$ and any sequence ${\{f_k\}_{k=1}^n}\subset X$ of
mean zero
i.r.v.'s we have the following sharp estimate:
%
%
\begin{equation}\label{XRosenthal}
\Biggl\|\sum_{k=1}^n f_k\Biggr\|_X\leq\alpha\|K\|_{X\to X}
\Biggl(\bigl\|F^*\chi_{[0,1]}\bigr\|_X+ \Biggl(\sum
_{k=1}^nF^*(k)^2 \Biggr)^{1/2} \Biggr).
\end{equation}
\end{proposition}
\end{remark}

The second main result of this section yields optimal (in order)
constant in a somewhat more general setting recently studied in \cite{M,Ju}.
\begin{theorem}\label{mainLptheorem} If $p,q\in[1,\infty)$, then for
an arbitrary sequence
$\{f_k\}_{k=1}^\infty\subset L_p$ of i.r.v.'s we have
%
%
\begin{equation}\label{mainLpestimate}\quad
\| \|\mathbf{f}\|_{l_q}\|_{L_p}\leq
\alpha\biggl(\frac{p} {\ln(p+1)} \biggr)^{1/q} \Biggl(\bigl\Vert
F^*\chi_{[0,1]}\bigr\Vert_{L_p} + \Biggl(\sum_{k=1}^\infty
F^*(k)^q \Biggr)^{1/q} \Biggr),
\end{equation}
where $\alpha$ is a universal
constant. Furthermore, the order $ (\frac{p}
{\ln(p+1)} )^{1/q}$ is optimal when $p\to\infty$.
\end{theorem}
\begin{pf} The fact that (\ref{mainLpestimate}) holds follows from a
combination of Theorem \ref{precise2} and Corollary
\ref{corollaryLp} above. It remains to show that the order
$ (\frac{p} {\ln(p+1)} )^{1/q}$ is optimal. Let
$f:=\chi_{[0,u]}$ $(0<u\leq1)$, $n\in\mathbb{N}$ and let
$\{f_{n,k}\}_{k=1}^n$ be a sequence of i.r.v.'s equidistributed
with the function $\sigma_{1/n}f=\chi_{[0,u/n]}$. Setting
$f_n:=\sum_{k=1}^nf_{n,k}$, we see that
\[
\lambda\{t\in[0,1]\dvtx f_n(t)=k\}=\frac{n!}{k!(n-k)!}\cdot
\biggl(\frac un \biggr)^k\cdot\biggl(1-\frac un \biggr)^{n-k},\qquad
k=1,2,\ldots,n.
\]
Hence,
\[
f_n^*(t)=\sum_{k=1}^n \chi_{[0, \tau_k^n]}\qquad \mbox{where }
\tau_k^n:=\sum_{i=k}^n \frac{n!}{i!(n-i)!}\cdot\biggl(\frac
un \biggr)^i\cdot\biggl(1-\frac un \biggr)^{n-i}.
\]
Note that for every $1\leq q<\infty$, we have
\[
\Biggl(\sum_{k=1}^nf^q_{n,k} \Biggr)^{1/q}= \Biggl(\sum
_{k=1}^nf_{n,k} \Biggr)^{1/q}=f_n^{1/q},
\]
and so the function $ (\sum_{k=1}^nf^q_{n,k} )^{1/q}$ is
equidistributed with the function
\[
(f_n^*)^{1/q}= \Biggl(\sum_{k=1}^n \chi_{[0,
\tau_k^n]} \Biggr)^{1/q}=\sum_{k=1}^n \bigl(k^{1/q}-(k-1)^{1/q}\bigr)\chi_{[0,
\tau_k^n]}.
\]
By the definition of the norm in Marcinkiewicz spaces, we have
%
%
\begin{eqnarray}\label{aux1}\qquad
\| \|(f_{n,k})\|_{l_q}\|_{L_p}
&\ge&\|(f_n^*)^{1/q}\|_{M(t^{1/p})}\ge\sup_{1\leq k\leq
n}\sum_{i=1}^k
\bigl(i^{1/q}-(i-1)^{1/q}\bigr)(\tau_k^n)^{1/p}\nonumber\\[-8pt]\\[-8pt]
&=&\sup_{1\leq k\leq n}k^{1/q}(\tau_k^n)^{1/p}.\nonumber
\end{eqnarray}
Estimating $\tau_k^n$ via Stirling's formula, we have
\begin{eqnarray*}
\tau_k^n &\ge&\frac{n!}{k!(n-k)!}\cdot\biggl(\frac un \biggr)^k\cdot
\biggl(1-\frac un \biggr)^{n-k}\\
&\ge&
\frac{\sqrt{n}n^ne^{-n}}{k!(n-k)^{n-k}e^{-n+k}\sqrt{n-k}}\cdot
\biggl(\frac un \biggr)^k\cdot\biggl(1-\frac un \biggr)^{n-k}\\
&\ge&
\frac{n^{n-k}u^ke^{-1}}{k! (n-k)^{n-k}e^k}\ge\frac
{1}{\sqrt{2\pi}e}\cdot\frac{u^k}{k^k\sqrt{k}}\ge\frac
{1}{\sqrt{2\pi}e}\cdot\frac{u^k}{k^{2k}}.
\end{eqnarray*}
Using the latter estimate in (\ref{aux1}), we obtain
%
%
\begin{equation}\label{aux2}
\| \|(f_{n,k})\|_{l_q}\|_{L_p}\ge\|(f_n^*)^{1/q}\|_{M(t^{1/p})}\ge
\frac{1}{\sqrt{2\pi}e}\sup_{1\leq k\leq n} k^{1/q} \biggl(\frac
{u^k}{k^{2k}} \biggr)^{1/p}.
\end{equation}

Observe further that for the sequence $\{f_{n,k}\}_{k=1}^n$ we have
$(\sum_{k=1}^n \bar f_{n,k})^*=\chi_{[0,u]}$ and therefore the term
$ (\Vert F^*\chi_{[0,1]}\Vert_{L_p}+ (\sum_{k=1}^\infty
F^*(k)^q )^{1/q} )$ in the right-hand side of
(\ref{mainLpestimate}) is equal in this case to
$\|\chi_{[0,u]}\|_{L_p}$. Thus, by (\ref{aux2}), we now estimate the
constant in (\ref{mainLpestimate}) from below as follows:
\begin{eqnarray*}
\sup_{0<u\leq1}\sup_{n\ge1}\frac{\|
\|(f_{n,k})\|_{l_q}\|_{L_p}}{\|\chi_{[0,u]}\|_{L_p}}&\ge&\sup
_{0<u\leq1}\sup_{n\ge1}\|(f_n^*)^{1/q}\|_{M(t^{1/p})}u^{-1/p}\\
&\ge& \frac{1}{\sqrt{2\pi}e}\sup_{0<u\leq1}\sup_{n\ge
1}n^{1/q} \biggl(\frac{u^n}{n^{2n}} \biggr)^{1/p} u^{-1/p}\\
&=&\frac
{1}{\sqrt{2\pi}e}\sup_{n\ge
1}\frac{n^{1/q}}{ ({n^{2n}} )^{1/p}}.
\end{eqnarray*}
%
Choosing $n= [\frac{p}{\ln p} ]$, we obtain for all
sufficiently large $p$'s
\[
\sup_{n\ge1}\frac{n^{1/q}}{ ({n^{2n}} )^{1/p}}\ge
\biggl(\frac{p}{\ln p}-1 \biggr)^{1/q}\cdot\biggl(\frac{\ln
p}{p} \biggr)^{{2}/{\ln(p)}}\ge\frac{1}{2e^{2}}
\biggl(\frac{p}{\ln(p+1)} \biggr)^{1/q}.
\]
The foregoing estimates show that the order $ (\frac{p}{\ln
(p+1)} )^{1/q}$ in (\ref{mainLpestimate}) is the best possible.
\end{pf}
\begin{remark}\label{mainLptheorem2} Since $\Lambda(t^{1/p})\subset
L_p\subset M(t^{1/p})$ ($1<p<\infty$),
it follows from (\ref{mainLpestimate}) that for an arbitrary
sequence $\{f_k\}_{k=1}^\infty\subset\Lambda(t^{1/p})$ of
i.r.v.'s
%
%
\begin{eqnarray}\label{mainLpestimate2}
\| \|\mathbf{f}\|_{l_q}\|_{M(t^{1/p})}&\leq&\alpha\biggl(\frac{p}
{\ln(p+1)} \biggr)^{1/q}\nonumber\\[-8pt]\\[-8pt]
&&{}\times \Biggl(\bigl\Vert F^*\chi_{[0,1]}\bigr\Vert
_{\Lambda(t^{1/p})}+ \Biggl(\sum_{k=1}^\infty
F^*(k)^q \Biggr)^{1/q} \Biggr),\nonumber
\end{eqnarray}
where $\alpha$ is a universal
constant. The argument used in the proof of the preceding
theorem shows that the order $ (\frac{p}
{\ln(p+1)} )^{1/q}$ remains optimal when $p\to\infty$.
\end{remark}

The following corollary from Theorem \ref{mainLptheorem} strengthens
\cite{AS1}, Theorem 4.4.
\begin{corollary}\label{newTh4.4} Let $1\leq q<\infty$. If a sequence
$\{f_k\}_{k\ge1}$ of uniformly bounded i.r.v.'s satisfies the
assumptions
%
\[
{\sup_{k\ge1}\sup_{t\in
[0,1]}}|f_k(t)|<\infty \quad\mbox{and}\quad \sum_{k=1}^\infty
F^*(k)^q<\infty,
\]
then the function $ ({\sum_{k=1}^\infty}|f_k|^q )^{1/q}\in
L_{N_q}$, where $L_{N_q}$ is the Orlicz space on $[0,1]$ generated
by the function $N_q(t):=t^{t^q}-1$, $N_q(0)=0$.
\end{corollary}
\begin{pf}
Due to (\ref{mainLpestimate}), we have
\[
\sup_{p\ge1} \biggl(\frac{\ln(p+1)}{p} \biggr)^{1/q}\| \|\mathbf{f}\|_{l_q}\|_{L_p}\leq\alpha
\Biggl(\bigl\|F^*\chi_{[0,1]}\bigr\|_\infty+ \Biggl(\sum_{k=1}^\infty
F^*(k)^q \Biggr)^{1/q} \Biggr),
\]
and, by the assumptions, the right-hand side of this inequality is
finite. Since the left-hand side coincides (up to equivalence) with
the norm of the function $\|\mathbf{f}\|_{l_q}= ({\sum_{k=1}^\infty
}|f_k|^q )^{1/q}$ in the Orlicz space $L_{N_q}$ (see, e.g.,
\cite{A}, Corollary 1, or~\cite{JSZ}, Proposition 3.6), and we are done.
\end{pf}

\section{Two complements to Rosenthal's inequality}\label{sec5}

The main results of this section are Proposition \ref{Sversion} and
Theorem \ref{Uversion} which complement Theorem \ref{precise1}.

Let ${\{f_k\}_{k=1}^n}$ be a sequence of i.r.v.'s on $[0,1]$,
$S_n:=\sum_{k=1}^n f_k$ and, as before, $F$ be the
disjointification function related to the sequence ${\{f_k\}_{k=1}^n}$ [see
(\ref{f1})], which may be written in the form
\[
F(t):=\sum_{k=1}^n f_k(t-k+1)\chi_{[k-1,k]}(t),\qquad t>0.
\]
\begin{proposition}\label{Sversion}
Let $X$ be a symmetric space on $[0,1]$ such that the Kruglov operator
$K$ acts boundedly on $X$. Then there exists
a universal constant $\beta>0$ such that for every $n\in\mathbb{N}$
and any sequence
${\{f_k\}_{k=1}^n}\subset X$ of i.r.v.'s the following inequality holds:
%
%
\begin{equation}\label{f11}
\|S_n\|_X\le\beta\|K\|_{X\to X}\bigl(\bigl\|F^*\chi_{[0,1]}\bigr\|_X+\|S_n\|_{L_1}\bigr).
\end{equation}
\end{proposition}
\begin{pf}
First, we assume that
%
%
\begin{equation}\label{f12}
\int_0^1 f_k(s) \,ds=0\qquad (k=1,2,\ldots,n).
\end{equation}
In this case, from Proposition \ref{meanzeroXcase} (see also
\cite{AS3}, Theorem 3.1) we infer that
\[
\|S_n\|_X\le
\gamma\|K\|_{X\to X}\bigl(\bigl\|F^*\chi_{[0,1]}\bigr\|_X+\bigl\|F^*\chi_{[1,\infty)}\bigr\|_{L_2}\bigr),
\]
where $\gamma$ is a universal constant. In addition, for the space
$L_1$, we have by \cite{JS}, Theorem 1,
\[
\|S_n\|_{L_1}\ge
c\bigl(\bigl\|F^*\chi_{[0,1]}\bigr\|_{L_1}+\bigl\|F^*\chi_{[1,\infty)}\bigr\|_{L_2}\bigr).
\]
Combining these two estimates we obtain (\ref{f11}) under assumption
(\ref{f12}).

Suppose now that the sequence
${\{f_k\}_{k=1}^n}\subset X$ is an arbitrary sequence of i.r.v.'s. Setting
%
%
\begin{equation}\label{f13}
g_k:=f_k-\int_0^1 f_k(s) \,ds\qquad (k=1,2,\ldots,n)
\end{equation}
we obtain a sequence ${\{g_k\}_{k=1}^n}\subset X$ satisfying (\ref
{f12}), and,
therefore, by the above
%
%
\begin{equation}\label{f14}
\|\sigma_n\|_X\le\gamma\|K\|_{X\to X}\bigl(\bigl\|G^*\chi_{[0,1]}\bigr\|_X+\|
\sigma_n\|_{L_1}\bigr),
\end{equation}
where
%
%
\begin{equation}\label{f15}
\sigma_n= \sum_{k=1}^n g_k,\qquad G(t):=\sum_{k=1}^n
g_k(t-k+1)\chi_{[k-1,k]}(t),\qquad t>0.
\end{equation}
On one hand,
\begin{eqnarray*}
G(t) &=& \sum_{k=1}^n \biggl(f_k(t-k+1)-\int_0^1
f_k(s) \,ds \biggr)\chi_{[k-1,k]}(t)\\
&=& F(t)-\sum_{k=1}^n\int_0^1 f_k(s) \,ds\chi_{[k-1,k]}(t),
\end{eqnarray*}
which implies
\[
G^*\chi_{[0,1]}(t)\le
F^*\chi_{[0,1]}(t)+\max_{k=1,2,\ldots,n}\|f_k\|_{L_1},
\]
and,
therefore, in view of the embedding $X\subset L_1$ with the constant
1,
%
%
\begin{equation}\label{f16}
\bigl\|G^*\chi_{[0,1]}\bigr\|_X\le
\bigl\|F^*\chi_{[0,1]}\bigr\|_X+{\max_{k=1,2,\ldots,n}}\|f_k\|_{X}\le
2\bigl\|F^*\chi_{[0,1]}\bigr\|_X.
\end{equation}
On the other hand,
\[
\|\sigma_n\|_X\ge\|S_n\|_X- \biggl|\int_0^1S_n(u) \,du \biggr|\ge
\|S_n\|_X-\|S_n\|_{L_1}
\]
and
\[
\|\sigma_n\|_{L_1}\le\|S_n\|_{L_1}+ \biggl|\int_0^1S_n(u) \,du \biggr|\le
2\|S_n\|_{L_1}.
\]
Combining these estimates with (\ref{f14}), we
obtain
\[
\|S_n\|_X \le(2\gamma\|K\|_{X\to X}+1)\bigl(\bigl\|F^*\chi_{[0,1]}\bigr\|_X+\|S_n\|_{L_1}\bigr),
\]
and the assertion is established in view of the fact $\|K\|_{X\to
X}\ge1/e$.
\end{pf}
\begin{remark} The converse inequality to (\ref{f11}) fails in general.
However, if in addition $f_k\ge0$, $k=1,2,\ldots, n$, then for any
symmetric space $X$
\[
c\bigl(\bigl\|F^*\chi_{[0,1]}\bigr\|_X+\|S_n\|_{L_1}\bigr)\leq\|S_n\|_X
\]
for some universal constant $c>0$.
\end{remark}

Denote
\[
U_n(t):={\max_{k=1,2,\ldots,n}}|S_k(t)|,\qquad t\in[0,1].
\]
\begin{theorem}\label{Uversion}
Let us assume that $X$ is an interpolation space for the couple $(L_1,
L_\infty)$ and that the Kruglov operator $K$ acts boundedly on $X$.
Then there exists a
universal constant $\alpha>0$ such that for all $n\in\mathbb{N}$
and any sequence
${\{f_k\}_{k=1}^n}\subset X$ of i.r.v.'s the following inequality holds
%
%
\begin{eqnarray}\label{f17}
&&\tfrac15
\bigl(\bigl\|F^*\chi_{[0,1]}\bigr\|_X+\|U_n\|_{L_1}\bigr)\nonumber\\[-8pt]\\[-8pt]
&&\qquad\leq\|U_n\|_X\leq
\alpha\|K\|_{X\to X} \bigl(\bigl\|F^*\chi_{[0,1]}\bigr\|_X+\|U_n\|_{L_1}
\bigr).\nonumber
\end{eqnarray}
\end{theorem}
\begin{pf}
First, it is obvious that
\[
2U_n(t)\ge M_n(t):={\max_{k=1,2,\ldots,n}}|f_k(t)|.
\]
Appealing to (\ref{Proposition 2.1HM}), we see that
\[
M_n^*(t/2)\ge F^*(t)\qquad (0<t\le1),
\]
and so
\[
\|U_n\|_X\ge\tfrac14\bigl\|F^*\chi_{[0,1]}\bigr\|_X,
\]
whence
\[
\|U_n\|_X\ge\tfrac15\bigl(\bigl\|F^*\chi_{[0,1]}\bigr\|_X+\|U_n\|_{L_1}\bigr).
\]

Let us prove the right-hand side inequality in (\ref{f17}). It holds
if $f_k$'s are symmetrically distributed. Indeed, by the well-known
Levy theorem (see, e.g., \cite{KW}, Proposition 1.1.1), we
have
\[
\lambda\{t\in[0,1]\dvtx U_n(t)>\tau\}\le2\lambda\{t\in
[0,1]\dvtx |S_n(t)|>\tau\} \qquad(\tau> 0),
\]
and so the result follows
from (\ref{f11}) (with the universal constant $\beta$) and the
assumption that $X$ is a symmetric space. Moreover, by the
standard ``symmetrization trick'' and using the assumption that $X$
is an interpolation space for the couple $(L_1, L_\infty)$, it is
not hard to extend this result to all sequences ${\{f_k\}
_{k=1}^n}\subset X$ of
i.r.v.'s satisfying condition (\ref{f12}).

Finally, let us consider the case of an arbitrary sequence
${\{f_k\}_{k=1}^n}\subset X$ of i.r.v.'s. Suppose that the sequence ${\{
g_k\}_{k=1}^n}$ is
defined by formula (\ref{f13}). Applying the (already established)
assertion to this sequence, we obtain
%
%
\begin{equation}\label{f20}
\|W_n\|_X\le2\beta\|K\|_{X\to X}\bigl(\bigl\|G^*\chi_{[0,1]}\bigr\|_X+\|W_n\|_{L_1}\bigr),
\end{equation}
where $G$ and $\sigma_k$ are defined as in (\ref{f15}), and
$W_n:=\max_{k=1,2,\ldots,n}|\sigma_k|$. Noting that
\[
W_n(t)\le U_n(t)+ {\max_{k=1,2,\ldots,n}\int_0^1}|S_k(u)| \,du\le
U_n(t)+\|U_n\|_{L_1},
\]
we infer $\|W_n\|_{L_1}\le2\|U_n\|_{L_1}$. Since
\begin{eqnarray*}
W_n(t) &=&
\max_{k=1,2,\ldots,n} \biggl|S_k(t)-\int_0^1S_k(u) \,du \biggr|\\
&\ge&
U_n(t) - {\int_0^1\max_{k=1,2,\ldots,n}}|S_k(u)| \,du=U_n(t)-\|U_n\|_{L_1},
\end{eqnarray*}
we have
\[
\|W_n\|_X\ge\|U_n\|_X-\|U_n\|_{L_1}.
\]
By (\ref{f20}) and
(\ref{f16}), this guarantees
\[
\|U_n\|_X\le(4\beta\|K\|_{X\to
X}+1)\bigl(\bigl\|F^*\chi_{[0,1]}\bigr\|_X+\|U_n\|_{L_1}\bigr),
\]
and the assertion
follows in view of the fact $\|K\|_{X\to X}\ge1/e$.
\end{pf}
\begin{remark} In the case when $X$ is a symmetric space containing $L_p$
for some finite $p$ (this condition is more restrictive than
the boundedness of the Kruglov operator in $X$, see
\cite{AS1}), the last result was also obtained in \cite{HM},
Theorem 5. At the same time, it is established in \cite{HM} even
for quasi-normed spaces.
\end{remark}

\section*{Acknowledgments}
The authors thank D. Zanin for a number of comments on an earlier
version of this paper. Also authors acknowledge support from the ARC.

%

%
\printaddresses


\begin{thebibliography}{23}

\bibitem{A}
%
\begin{barticle}[mr]
\bauthor{\bsnm{Astashkin},~\bfnm{S.~V.}\binits{S.~V.}}
(\byear{2005}).
\btitle{Extrapolation functors on a family of scales generated by the real
interpolation method}.
\bjournal{Sibirsk. Mat. Zh.}
\bvolume{46}
\bpages{264--289}.
\bid{doi={10.1007/s11202-005-0021-2}, mr={2141194}}
\end{barticle}
%
\endbibitem

\bibitem{AS1}
%
\begin{barticle}[mr]
\bauthor{\bsnm{Astashkin},~\bfnm{S.~V.}\binits{S.~V.}} \AND
\bauthor{\bsnm{Sukochev},~\bfnm{F.~A.}\binits{F.~A.}}
(\byear{2005}).
\btitle{Series of independent random variables in rearrangement invariant
spaces: An operator approach}.
\bjournal{Israel J. Math.}
\bvolume{145}
\bpages{125--156}.
\bid{doi={10.1007/BF02786688}, mr={2154724}}
\end{barticle}
%
\endbibitem

\bibitem{AS2}
%
\begin{barticle}[mr]
\bauthor{\bsnm{Astashkin},~\bfnm{S.~V.}\binits{S.~V.}} \AND
\bauthor{\bsnm{Sukochev},~\bfnm{F.~A.}\binits{F.~A.}}
(\byear{2004}).
\btitle{Comparison of sums of independent and disjoint functions in symmetric
spaces}.
\bjournal{Mat. Zametki}
\bvolume{76}
\bpages{483--489}.
\bid{doi={10.1023/B:MATN.0000043474.00734.ec}, mr={2112064}}
\end{barticle}
%
\endbibitem

\bibitem{AS3}
%
\begin{barticle}[vtex]
\bauthor{\bsnm{Astashkin},~\bfnm{S.~V.}\binits{S.~V.}} \AND
\bauthor{\bsnm{Sukochev},~\bfnm{F.~A.}\binits{F.~A.}}
(\byear{2008}).
\btitle{Series of independent mean zero random variables in
rearrangement-invariant spaces having the Kruglov
property}.
\bjournal{J. Math. Sci.}
\bvolume{148}
\bpages{795--809}.
\end{barticle}
%
\endbibitem

\bibitem{Br}
%
\begin{bbook}[mr]
\bauthor{\bsnm{Braverman},~\bfnm{Michael~Sh.}\binits{M.~S.}}
(\byear{1994}).
\btitle{Independent Random Variables and Rearrangement Invariant Spaces}.
\bseries{London Mathematical Society Lecture Note Series}
\bvolume{194}.
\bpublisher{Cambridge Univ. Press}, \baddress{Cambridge}.
\bid{mr={1303591}}
\end{bbook}
%
\endbibitem

\bibitem{Bu}
%
\begin{barticle}[mr]
\bauthor{\bsnm{Bukhvalov},~\bfnm{A.~V.}\binits{A.~V.}}
(\byear{1987}).
\btitle{Interpolation of linear operators in spaces of vector
functions and
with a mixed norm}.
\bjournal{Sibirsk. Mat. Zh.}
\bvolume{28}
\bpages{37--51}.
\bid{mr={886851}}
\end{barticle}
%
\endbibitem

\bibitem{Bur}
%
\begin{barticle}[mr]
\bauthor{\bsnm{Burkholder},~\bfnm{D.~L.}\binits{D.~L.}}
(\byear{1979}).
\btitle{A sharp inequality for martingale transforms}.
\bjournal{Ann. Probab.}
\bvolume{7}
\bpages{858--863}.
\bid{mr={542135}}
\end{barticle}
%
\endbibitem

\bibitem{Cal}
%
\begin{barticle}[mr]
\bauthor{\bsnm{Calder{\'o}n},~\bfnm{A.~P.}\binits{A.~P.}}
(\byear{1964}).
\btitle{Intermediate spaces and interpolation, the complex method}.
\bjournal{Studia Math.}
\bvolume{24}
\bpages{113--190}.
\bid{mr={0167830}}
\end{barticle}
%
\endbibitem

\bibitem{Cr}
%
\begin{barticle}[mr]
\bauthor{\bsnm{Creekmore},~\bfnm{J.}\binits{J.}}
(\byear{1981}).
\btitle{Type and cotype in {L}orentz {$L\sb{pq}$} spaces}.
\bjournal{Nederl. Akad. Wetensch. Indag. Math.}
\bvolume{43}
\bpages{145--152}.
\bid{mr={707247}}
\end{barticle}
%
\endbibitem

\bibitem{HM}
%
\begin{barticle}[mr]
\bauthor{\bsnm{Hitczenko},~\bfnm{Pawe{\l}}\binits{P.}} \AND
\bauthor{\bsnm{Montgomery-Smith},~\bfnm{Stephen}\binits{S.}}
(\byear{2001}).
\btitle{Measuring the magnitude of sums of independent random variables}.
\bjournal{Ann. Probab.}
\bvolume{29}
\bpages{447--466}.
\bid{doi={10.1214/aop/1008956339}, mr={1825159}}
\end{barticle}
%
\endbibitem

\bibitem{JMST}
%
\begin{barticle}[mr]
\bauthor{\bsnm{Johnson},~\bfnm{W.~B.}\binits{W.~B.}},
\bauthor{\bsnm{Maurey},~\bfnm{B.}\binits{B.}},
\bauthor{\bsnm{Schechtman},~\bfnm{G.}\binits{G.}} \AND
\bauthor{\bsnm{Tzafriri},~\bfnm{L.}\binits{L.}}
(\byear{1979}).
\btitle{Symmetric structures in {B}anach spaces}.
\bjournal{Mem. Amer. Math. Soc.}
\bvolume{19}
\bpages{v+298}.
\bid{mr={527010}}
\end{barticle}
%
\endbibitem

\bibitem{JS}
%
\begin{barticle}[mr]
\bauthor{\bsnm{Johnson},~\bfnm{William~B.}\binits{W.~B.}} \AND
\bauthor{\bsnm{Schechtman},~\bfnm{G.}\binits{G.}}
(\byear{1989}).
\btitle{Sums of independent random variables in rearrangement invariant
function spaces}.
\bjournal{Ann. Probab.}
\bvolume{17}
\bpages{789--808}.
\bid{mr={985390}}
\end{barticle}
%
\endbibitem

\bibitem{JSZ}
%
\begin{barticle}[mr]
\bauthor{\bsnm{Johnson},~\bfnm{W.~B.}\binits{W.~B.}},
\bauthor{\bsnm{Schechtman},~\bfnm{G.}\binits{G.}} \AND
\bauthor{\bsnm{Zinn},~\bfnm{J.}\binits{J.}}
(\byear{1985}).
\btitle{Best constants in moment inequalities for linear combinations of
independent and exchangeable random variables}.
\bjournal{Ann. Probab.}
\bvolume{13}
\bpages{234--253}.
\bid{mr={770640}}
\end{barticle}
%
\endbibitem

\bibitem{Ju}
%
\begin{barticle}[mr]
\bauthor{\bsnm{Junge},~\bfnm{Marius}\binits{M.}}
(\byear{2006}).
\btitle{The optimal order for the {$p$}-th moment of sums of
independent random
variables with respect to symmetric norms and related combinatorial
estimates}.
\bjournal{Positivity}
\bvolume{10}
\bpages{201--230}.
\bid{doi={10.1007/s11117-005-0008-z}, mr={2237498}}
\end{barticle}
%
\endbibitem

\bibitem{K}
%
\begin{barticle}[mr]
\bauthor{\bsnm{Kruglov},~\bfnm{V.~M.}\binits{V.~M.}}
(\byear{1970}).
\btitle{A remark on the theory of infinitely divisible laws}.
\bjournal{Teor. Veroyatn. Primen.}
\bvolume{15}
\bpages{330--336}.
\bid{mr={0272018}}
\end{barticle}
%
\endbibitem

\bibitem{KPS}
%
\begin{bbook}[mr]
\bauthor{\bsnm{Kre{\u\i}n},~\bfnm{S.~G.}\binits{S.~G.}},
\bauthor{\bsnm{Petun{\={\i}}n},~\bfnm{Yu.~{\=I}.}\binits{Y.~{\=
I}.}} \AND
\bauthor{\bsnm{Sem{\"e}nov},~\bfnm{E.~M.}\binits{E.~M.}}
(\byear{1982}).
\btitle{Interpolation of Linear Operators}.
\bseries{Translations of Mathematical Monographs}
\bvolume{54}.
\bpublisher{Amer. Math. Soc.}, \baddress{Providence, RI.}
\bid{mr={649411}}
\end{bbook}
%
\endbibitem

\bibitem{KW}
%
\begin{bbook}[mr]
\bauthor{\bsnm{Kwapie{\'n}},~\bfnm{Stanis{\l}aw}\binits{S.}} \AND
\bauthor{\bsnm{Woyczy{\'n}ski},~\bfnm{Wojbor~A.}\binits{W.~A.}}
(\byear{1992}).
\btitle{Random Series and Stochastic Integrals: Single and Multiple}.
\bpublisher{Birkh\"auser}, \baddress{Boston, MA}.
\bid{mr={1167198}}
\end{bbook}
%
\endbibitem

\bibitem{La}
%
\begin{barticle}[mr]
\bauthor{\bsnm{Lata{\l}a},~\bfnm{Rafa{\l}}\binits{R.}}
(\byear{1997}).
\btitle{Estimation of moments of sums of independent real random variables}.
\bjournal{Ann. Probab.}
\bvolume{25}
\bpages{1502--1513}.
\bid{doi={10.1214/aop/1024404522}, mr={1457628}}
\end{barticle}
%
\endbibitem

\bibitem{Lo}
%
\begin{barticle}[mr]
\bauthor{\bsnm{Lozanovski{\u\i}},~\bfnm{G.~Ja.}\binits{G.~J.}}
(\byear{1972}).
\btitle{A remark on a certain interpolation theorem of {C}alder\'on}.
\bjournal{Funktsional. Anal. i Prilozhen.}
\bvolume{6}
\bpages{89--90}.
\bid{mr={0312246}}
\end{barticle}
%
\endbibitem

\bibitem{LT2}
%
\begin{bbook}[mr]
\bauthor{\bsnm{Lindenstrauss},~\bfnm{Joram}\binits{J.}} \AND
\bauthor{\bsnm{Tzafriri},~\bfnm{Lior}\binits{L.}}
(\byear{1979}).
\btitle{Classical {B}anach Spaces. {II}. Function Spaces}.
\bseries{Ergebnisse der Mathematik und Ihrer Grenzgebiete [Results in
Mathematics and Related Areas]}
\bvolume{97}.
\bpublisher{Springer}, \baddress{Berlin}.
\bid{mr={540367}}
\end{bbook}
%
\endbibitem

\bibitem{M}
%
\begin{barticle}[mr]
\bauthor{\bsnm{Montgomery-Smith},~\bfnm{Stephen}\binits{S.}}
(\byear{2002}).
\btitle{Rearrangement invariant norms of symmetric sequence norms of
independent sequences of random variables}.
\bjournal{Israel J. Math.}
\bvolume{131}
\bpages{51--60}.
\bid{doi={10.1007/BF02785850}, mr={1942301}}
\end{barticle}
%
\endbibitem

\bibitem{Pr}
%
\begin{barticle}[mr]
\bauthor{\bsnm{Prohorov},~\bfnm{Yu.~V.}\binits{Y.~V.}}
(\byear{1958}).
\btitle{Strong stability of sums and infinitely divisible laws}.
\bjournal{Teor. Veroyatn. Primen.}
\bvolume{3}
\bpages{153--165}.
\bid{mr={0098428}}
\end{barticle}
%
\endbibitem

\bibitem{Ros}
%
\begin{barticle}[mr]
\bauthor{\bsnm{Rosenthal},~\bfnm{Haskell~P.}\binits{H.~P.}}
(\byear{1970}).
\btitle{On the subspaces of {$L_{p}$} {$(p>2)$} spanned by
sequences of
independent random variables}.
\bjournal{Israel J. Math.}
\bvolume{8}
\bpages{273--303}.
\bid{mr={0271721}}
\end{barticle}
%
\endbibitem

\end{thebibliography}
\end{document}